\documentclass[12pt,a4paper]{article}
\usepackage{amssymb}

\usepackage{graphicx,amssymb,amsfonts,epsfig,amsthm,a4,amsmath,url}
\usepackage[latin1]{inputenc}

\psfigdriver{dvips}

\newtheorem{thm}{Theorem}
\newtheorem{cor}[thm]{Corollary}
\newtheorem{lem}[thm]{Lemma}

\newtheorem{prop}[thm]{Proposition}

\theoremstyle{definition}
\newtheorem{defn}[thm]{Definition}

\theoremstyle{remark}
\newtheorem{rem}[thm]{Remark}
\newtheorem{exe}[thm]{Example}
\newtheorem{exes}[thm]{Examples}

\numberwithin{equation}{section}

\newcommand{\A}{\mathcal{A}}

\newcommand{\HH}{\mathcal{H}}

\newcommand{\N}{\mathbf{N}}

\newcommand{\U}{\mathcal{U}}
\newcommand{\F}{\mathcal{F}}
\newcommand{\M}{\mathcal{M}}
\newcommand{\R}{\mathbf{R}}

\newcommand{\V}{\mathbf{Var}}
\newcommand{\E}{\mathbf{E}}

\newcommand{\tpr}{\begin{tiny}\noindent Proof:}

\newcommand{\card}{\textnormal{card}}

\newcommand{\bpr}{\noindent \textbf{Proof}: ~}
\newcommand{\epr}{~$\blacksquare$}

\title{Coarse embeddings into a Hilbert space, Haagerup Property and Poincar\'e inequalities}%
\author{Romain Tessera\footnote{The author is supported by the
NSF grant DMS-0706486.}}
\date{\today}

\begin{document}

\baselineskip=16pt

\maketitle

 \maketitle
\begin{abstract}
We prove that a metric space does not coarsely embed into a Hilbert
space if and only if it satisfies a sequence of Poincar\'e
inequalities, which can be formulated in terms of (generalized)
expanders. We also give quantitative statements, relative to the
compression. In the equivariant context, our result says that a
group does not have the Haagerup property if and only if it has
relative property T with respect to a family of probabilities whose
supports go to infinity. We give versions of this result both in
terms of unitary representations, and in terms of affine isometric
actions on Hilbert spaces.
\end{abstract}

\section{Introduction}
\subsection{Obstruction to coarse embeddings}
The notion of expanders has been pointed out by Gromov as an
obstruction for a metric space to coarsely embed into a Hilbert
space. Recall \cite{JS} (see also \cite{Lub}) that a sequence of
expanders is a sequence of finite connected graphs $(X_n)$ with
bounded degree, satisfying the following Poincar\'e inequality for
all $f\in \ell^2(X_n)$

\begin{equation}\label{trueexpander}
    \frac{1}{|X_n|^2}\sum_{x,y\in X_n}|f(x)-f(y)|^2\leq \frac{C}{|X_n|}\sum_{x\sim y}|f(x)-f(y)|^2,
\end{equation}
for some constant $C>0$, and whose cardinality $|X_n|$ goes to
infinity when $n\to \infty$. An equivalent formulation in $\ell^p$
\cite{M2} can be used to prove that expanders do not coarsely embed
into $L^p$ for any $1\leq p<\infty.$

It is an open problem whether a metric space with bounded geometry
that does not coarsely embed into a Hilbert space admits a coarsely
embedded sequence of expanders.

In this paper, we prove that a metric space (not necessarily with
bounded geometry) that does not coarsely embed into a Hilbert space
admits a coarsely embedded sequence of ``generalized expanders".
This weaker notion of expanders can be roughly described as a
sequence of Poincar\'e inequalities with respect to finitely
supported probability measures on $X\times X$. We also provide
similar obstructions for coarse embeddings into families of metric
spaces such as $L^p$, for every $1\leq p<\infty$, uniformly convex
Banach spaces, and CAT(0) spaces.

For the sake of clarity, we chose to present most of our results
first in the case of Hilbert spaces. However, our characterization
(see Theorem~\ref{MainThm}) of the non-existence of coarse embedding
into $L^p$ deserves some attention. Indeed, our Poincar\'e
inequalities are not equivalent for different values of $2\leq p
<\infty$. This follows from a result of Naor and Mendel \cite{MN}
(see also \cite{JR}) saying that $L^p$ does not coarsely embed into
$L^q$ if $2\leq q<p$. This is different from what happens with real
expanders, as having a sequence of expanders prevents from having a
coarse embedding into $L^p$, for any $1\leq p<\infty$. In
particular, at least without any assumption of bounded geometry, our
generalized expanders cannot be replaced by actual expanders.

To conclude, let us remark that finding subspaces of $L^p$ for some
$p>2$, with bounded geometry, which do not coarsely embed into
$L^2$, would answer negatively the problem mentioned above.

\subsection{Obstruction to Haagerup Property}

A countable group is said to have the Haagerup property if it admits
a proper affine action on a Hilbert space. An obstruction for an
infinite countable group to have the Haagerup Property is known as
Property T (also called Property FH), which says that every
isometric affine action has a fixed point (or equivalently bounded
orbits). A weaker obstruction is to have relative property T with
respect to an infinite subset \cite{C,C'}. The case where this
subset is a normal subgroup has been mostly considered, as it has
strong consequences. On the other hand, there are examples of groups
which do not have relative property T with respect to any subgroup,
but have it with respect to some infinite subset \cite{C}. The
question whether the latter property is equivalent to the negation
of Haagerup Property is still open.

In this paper, we partially answer this question by showing that a
countable group does not have Haagerup Property if and only if it
has relative Property T with respect to a sequence of probabilities
whose supports eventually leave every finite subset.

\bigskip

\noindent \textbf{Acknowledgments.} I am grateful to Yann Ollivier,
Yves de Cornulier, and Bogdan Nica for their useful comments and
corrections. I address a special thank to James Lee who pointed to
me  \cite[Proposition~15.5.2]{M}.

\section{Statement of results in the Hilbert case}
In this section, we state our main results concerning embeddings
into a Hilbert space. In Section~\ref{MoregeneralresultsSection},
using a slightly more sophisticated vocabulary, we generalize to
other geometries.

\subsection{Coarse embeddability into Hilbert spaces and generalized expanders}
Let $\HH$ denote a separable infinite dimensional Hilbert space. We
denote by $|v|$ the norm of a vector in $\HH.$ Let $X=(X,d)$ be a
metric space. For all $r\geq 0$, denote $$\Delta_r(X)=\{(x,y)\in
X^2, d(x,y)\geq r\}.$$ In this paper, we prove that a metric space
that does not coarsely embed into a Hilbert space contains in a weak
sense a sequence of expanders. Precisely, following the idea  of
\cite[Section~4.2]{T}, let us define
\begin{defn}

\

\begin{itemize}
\item
Let $K$ and $r$ be positive numbers. A finite metric space is called
a {\it generalized $(K,r)$-expander} if there exists a symmetric
probability measure $\mu$ supported on $\Delta_{r}(X)$ with the
following property. For every map $F:X\to \HH$ satisfying
$|F(x)-F(y)|\leq d(x,y)$ for all $(x,y)\in \Delta_1(X)$, we have
    \begin{equation}\label{expander'}
\V_{\mu}(F):=\sum_{x,y}|F(x)-F(y)|^2\mu(x,y)\leq K^2.
    \end{equation}

\item A sequence of finite metric spaces $(X_n)$ is called a {\it sequence of generalized $K$-expander} if for every $n\in \N$, $X_n$ is a $(K,r_n)$-expanders, where $r_n\to \infty.$
\end{itemize}
\end{defn}

Recall that a family of metric spaces $(X_i)_{i\in I}$ coarsely
embeds into a metric space $Y$ if there exists a family $(F_i)$ of
uniformly coarse embeddings of $X_i$ into $Y$, i.e. if there are two
increasing, unbounded functions $\rho_-$ and $\rho_+$ such that
\begin{equation}\label{coarseEq}
\rho_-(d(x,y))\leq d(F_i(x),F_i(y))\leq \rho_+(d(x,y)), \forall
x,y\in X_i, \forall i\in I.
\end{equation}
\begin{prop}
A sequence of generalized expanders $(X_n)$ does not coarsely embed
into a Hilbert space.
\end{prop}\label{coearse_embedding_of_expanderprop}
\bpr Let $K>0$ and for all $n\in \N$, let $X_n$ is a
$(K,r_n)$-expander, with $r_n\to \infty$. For every $n\in \N$, let
$F_n$ be a map from $X_n\to \HH$, and that there exists some
increasing function $\rho_+$ such that $|F_n(x)-F_n(y)|\leq
\rho_+(d(x,y)), \forall x,y\in X_n, \forall n\in \N.$ As observed in
\cite[Lemmas~2.4 and~3.11]{CTV}, if a metric space (or a family of
metric spaces) coarsely embeds into a Hilbert space, we can always
assume that the function $\rho_+$ goes arbitrarily slowly to
infinity (this follows from a result of Bochner and Schoenberg
\cite[Theorem 8]{Sch}). So in particular, we can assume that
$\rho_+(t)\leq t, \; \forall t\geq 1$. But then, (\ref{expander'})
tells us that pairs of points of $X_n$, which are at distance $\geq
r_n$ are sent by $F_n$ at distance less than $K$. As $r_n\to
\infty$, this implies that any increasing function $\rho_-$
satisfying $$\rho_-(d(x,y))\leq |F_n(x)-F_n(y)|,\;\forall x,y\in
X_n, \forall n\in \N$$ would have to be $\leq K.$ \epr

\

Our main result is the following theorem (which is a particular case
of Corollary~\ref{Maincor}).

\begin{thm}
A metric space does not coarsely embed into a Hilbert space if and
only if it has a coarsely-embedded sequence of generalized
expanders.
\end{thm}

\subsection{Comparison with the usual notion of expanders}

The usual definition of an expander is a sequence of finite
connected graphs $(X_n)$  with degree $\leq k$, satisfying the
following Poincar\'e inequality for all $f\in \ell^2(X_n)$

\begin{equation}\label{trueexpander}
    \frac{1}{|X_n|^2}\sum_{x,y\in X_n}|f(x)-f(y)|^2\leq \frac{C}{|X_n|}\sum_{x\sim y}|f(x)-f(y)|^2,
\end{equation}
for some constant $C>0$, and whose cardinality $|X_n|$ goes to
infinity when $n\to \infty$. If $\nu_n$ denote the uniform measure
on $X_n\times X_n$, this can be rewritten as
$$\V_{\nu_n}(f)\leq \frac{C}{|X_n|}\sum_{x\sim y}|f(x)-f(y)|^2.$$
Now, assuming that $f$ is $1$-Lipschitz, we have
$$\V_{\nu_n}(f)\leq kC,$$
To obtain condition (\ref{expander'}), we need to replace $\nu_n$ by
a probability supported far away from the diagonal. To do that, we
just notice that at least half of the mass of $\nu_n$ is actually
supported on $\Delta_{r_n}$, with $r_n=\log_k(|X_n|/2)$. Indeed, if
$r$ is some positive number, the number of pairs of $X_n\times X_n$
which are at distance $\leq r$ is at most $k^r|X_n|$. Hence the
proportion of such pairs is $\leq k^r/|X_n|$, and the statement
follows. Therefore, renormalizing the restriction of $\nu_n$ to
$\Delta_{r_n}$, we obtain a probability $\mu_n$ satisfying
$\V_{\mu_n}(f)\leq 2\V_{\nu_n}(f)\leq 2kC.$ Hence, we have proved
\begin{prop}
A sequence of expanders satisfying (\ref{trueexpander}) with
constant $C$, is a sequence of generalized $K$-expanders, with
$K=(2kC)^{1/2}.$ \epr
\end{prop}

\subsection{Haagerup property and relative property T with respect to a family of probabilities}

Recall that a countable group has the Haagerup Property if it acts
metrically properly by affine isometries on a Hilbert space. It is
well known that this is equivalent to saying that $G$ has a proper
Hilbert length (see Section~\ref{MoregeneralresultsSection} for a
definition). On the other hand \cite{C} a group has relative
Property FH with respect to an infinite subset $\Omega$ if every
Hilbert length on $G$ is bounded in restriction to $\Omega$.

\begin{defn}
Let $G$ be a countable group equipped with a proper length function
$L_0$. Let $(\mu_n)$ be a sequence of probability measures on $G$.
We say that $G$ has {\it relative property FH with respect to
$(\mu_n)$} there exists $K>0$ such that for every Hilbert length $L$
satisfying
$$L(g)\leq L_0(g), \quad \forall g\in G,$$ and for every $n\in \N$,
$$\E_{\mu_n}(L^2)\leq K.$$
\end{defn}
It is easy to see that this definition does not depend on $L_0$.

Note that having relative property FH with respect to an infinite
subset $\Omega=\{a_1,a_2,\ldots\}$ corresponds to having relative
Property FH with respect to $(\mu_n)$, where $\mu_n$ is the Dirac
measure at $a_n$, for every $n\in \N$.

\begin{thm}\label{Haagerup}
A countable group $G$ does not have the Haagerup Property if and
only if it has relative Property FH with respect to a sequence of
symmetric probability measures $(\mu_n)$, such that for all $n\in
\N$, $\mu_n$ is supported on a finite subset of $\{g,L_0(g)\geq
n\}$.
\end{thm}

Recall that an equivalent formulation of the Haagerup Property
(actually the original one) is as follows: there exists a sequence
$(\phi_k)$ of positive definite functions on the group such that
$\lim_{k\to \infty}\phi_k(g)=1$ for all $g\in G$, and $\lim_{g\to
\infty}\phi_k(g)=0$ for all $k\in \N$ (in terms of unitary
representations, it says that there exists a $C_0$ unitary
representation with almost-invariant vectors).

An obvious obstruction to the Haagerup Property is \cite{C} relative
property T with respect to an unbounded subset $\Omega$: every
sequence $(\phi_k)$ of positive definite function on $G$ converging
to $1$ pointwise, converges uniformly in restriction to $\Omega$. In
\cite{C}, it is actually proved that relative Property T with
respect to $\Omega$ is equivalent to relative Property FH with
respect to $\Omega.$ Let us introduce the following definition.
\begin{defn}
Let $G$ be a countable group. Let $(\mu_n)$ be a sequence of
probability measures on $G$. We say that $G$ has {\it relative
property T with respect to $(\mu_n)$} if every sequence of positive
definite function $(\phi_k)$ that pointwise converges to $1$,
satisfies that $\lim_{k\to \infty}\E_{\mu_n}(\phi_k)=1$ uniformly
with respect to $n\in \N$.
\end{defn}

We have the following theorem

\begin{thm}\label{HaagerupBis}
A countable group $G$ does not have the Haagerup Property if and
only if it has relative Property T with respect to a sequence of
symmetric probability measures $(\mu_n)$, such that for all $n\in
\N$, $\mu_n$ is supported on a finite subset of $\{g,L_0(g)\geq
n\}$.
\end{thm}
\bpr It is clear that relative Property T with respect to a sequence
of probabilities whose supports go to infinity violates the Haagerup
Property. So what we need to prove is the converse, namely, that the
negation of Haagerup Property implies relative Property T with
respect to some $(\mu_n)$. By Theorem~\ref{Haagerup}, it is enough
to prove that relative Property FH implies relative Property T,
which is a straightforward adaptation of the proof of
\cite[Theorem~3]{AW}. \epr

\section{A more general setting, and quantitative statements}\label{MoregeneralresultsSection}

In this section, we switch to a slightly different point of view.
The statements we want to prove are of the following form: a metric
space $X$ that cannot coarsely embed into some class of spaces $\M$
has to satisfy some sequence of Poincar\'e inequalities. It is worth
noting that these inequalities consist essentially in a comparison
between metrics on $X$. Namely, we compare the original metric on
$X$ with all the pull-back metrics obtained from maps to metric
spaces of $\M$. let us be more precise.

Let $X$ be a set. A pseudo-metric on $X$ is a function: $\sigma:
X^2\to \R_+$ such that $\sigma(x,y)=\sigma(y,x)$, $\sigma(x,y)\leq
\sigma(x,z)+\sigma(z,y)$, and $\sigma(x,x)=0,$ for all $x,y,z\in X.$
In the sequel, a pseudo-metric will simply be called a metric.

If $(Y,d)$ is a metric space and $F:X\to Y$ is a map, then we can
consider the pull-back metric $\sigma_F(x,y)=d(F(x),F(y))$, for all
$x,y\in X$. Such metrics are called $Y$-metrics on $X$. More
generally, if $\M$ is a class of metric spaces, a $\M$-metric on $X$
is a $Y$-metric for some $Y\in \M.$

Assume here that $X=(X,d)$ is a metric space. A metric $\sigma$ is a
called coarse if there exist two increasing unbounded functions
$\rho_-, \rho_+$ such that, for all $x,y\in X,$ $$\rho_-(d(x,y))\leq
\sigma(x,y)\leq \rho_+(d(x,y)).$$ Note that if $\sigma=\sigma_F$ is
a $Y$-metric associated to a map $F:X\to Y$, then $\sigma_F$ is
coarse if and only if $F$ is a coarse embedding.
\begin{defn}
A {\it sheaf of metrics} on a set $X$ is a collection of pairs
$(\sigma,\Omega)$, where $\Omega$ is a subset of $X$, and $\sigma$
is a metric defined on $\Omega$. If $\Omega$ is an subset of $X$, we
denote by $\F(\Omega)$ the set of pairs $(\sigma,\Omega)\in \F$.

We also assume that the restriction is well-defined from
$\F(\Omega)$ to $\F(\Omega')$ for every $\Omega'\subset \Omega$
(which is automatic in the case of sheaves of $\M$-metrics).
\end{defn}

One checks easily that squares of Hilbert metrics, and more
generally $p$-powers of $L^p$-metrics form a convex sub-cone of the
space of real-valued functions on $X^2$. This is in fact a crucial
remark for what follows.

\begin{defn}
Let $X$ be a set. A sheaf $\F$ of metrics on $X$ is called {\it
$p$-admissible} (for some $p>0$) if for every $\Omega$, the
following hold.
\begin{itemize}\label{p-admissibleDef}
\item[(i)]
The set of $\sigma^p$, where $\sigma\in \F(\Omega)$ forms a sub-cone
of the space of functions on $\Omega^2$.

\item[(ii)]  $\F(\Omega)$ is closed for the topology of pointwise
convergence.

\item[(iii)]
Let $(U_i)$ be a family of finite subsets whose union is $\Omega$,
satisfying that for all $i,j\in I$, there exists $k$ such that
$U_i\cup U_j\subset U_k.$ Let $(\sigma_i,U_i)$ be a compatible
family of sections, in the sense that $\sigma_i$ and $\sigma_j$
coincide on the intersection $U_i\cap U_j$. Then there exists a
section $\sigma\in \F(\Omega)$, whose restriction to every $U_i$ is
$\sigma_i$. In other words, $\F(\Omega)$ is the direct limit of the
$\F(U_i)$.
\end{itemize}
\end{defn}

\begin{prop}
Let $X$ be a set, and let $\M$ be a class of metric spaces which is
closed under ultra-limits. Then the sheaf of $\M$-metrics on
(subsets of) $X$ satisfies conditions (ii) and (iii) of
Definition~\ref{p-admissibleDef}.
\end{prop}
\bpr That $X$ satisfies (ii) is trivial. Let us prove (iii). Fix a
point $o\in \Omega$ and consider $\A$ the partially ordered (for the
inclusion) set of all finite subsets of $\Omega$ containing $o$. For
every $\sigma_U\in \F(U)$, where $U\in \A$, choose $Y_F\in \M$ and
$F_U: U\to Y_U$ be such that $\sigma_U(x,y)=d(F_U(x),F_U(y))$ for
all $x,y\in U$. In every $Y_U$, take $y_U=F_U(o)$ for the origin.
Fix an ultra-filter $\U$ on $\A$. Now, the limit $F$ of the $F_U$ is
well defined from $$\Omega=\bigcup_{U\in A_o}U\to
\lim_{\U}(Y_U,y_U),$$ and $\sigma(x,y)=d(F(x),F(y))$ satisfies the
third condition of Definition~\ref{p-admissibleDef}. \epr

As a consequence of the proposition, we get the following examples.
\begin{exes}

\

\begin{itemize}

\item for $p\geq 1$, the sheaf of $L^p$-metrics is $p$-admissible \cite{H}.

\item Let $c>0$ and  $1<p<\infty$. The class $\M_{c,p}$ of $(c,p)$-uniformly convex Banach spaces, is the class of uniformly
convex Banach spaces whose moduli of convexity satisfy
$\delta(t)\geq ct^p$. The sheaf of $\M_{c,p}$-metrics is
$p$-admissible.

\item The sheaf of CAT(0)-metrics is $2$-admissible.

\end{itemize}
\end{exes}

If $(X,\mu)$ is a probability space and $f$ is a integrable function on $X$, we denote $\E_{\mu}(f)=\int f(x)d\mu(x)$.

\begin{thm}\label{MainThm}
Let $X$ be a metric space. Let $\F$ be a $p$-admissible sheaf of
metrics on $X$. Then there exists no coarse metric in $\F(X)$ if and
only if for every function $\rho_+:\R_+\to \R_+$, there exist $K>0,$
and a sequence of symmetric probability measures $(\mu_n)$ with the
following properties
\begin{itemize}
\item[-] for every
$n\in \N$, $\mu_n$ is supported on a finite subset $A_n$ of
$\Delta_n(X)$; \item[-] for every $n\in \N$ and every $\sigma\in
\F(A_n)$ satisfying
$$\sigma(x,y)\leq \rho_+(d(x,y)), \quad \forall (x,y)\in A_n,$$ one
has
\begin{equation}\label{maineq}
\E_{\mu_n}(\sigma^p)\leq K^p.
\end{equation}
\end{itemize}
\end{thm}

\begin{rem}
Note that this theorem characterizes metric spaces that do not
coarsely embed into $L^p$-spaces, $CAT(0)$-spaces, uniformly convex
spaces... Indeed, by a theorem of Pisier \cite{Pisier}, any
uniformly convex Banach space is isomorphic to a $(c,p)$-uniformly
convex Banach space for some $1<p<\infty$ and $c>0$. We can also
avoid to use this deep theorem by defining $\phi$-admissible sheafs
of metrics for any non-decreasing convex function $\phi$, and by
adapting the proof of Theorem~\ref{MainThm} to this slightly more
general setting.
\end{rem}

Generalizing the case of Hilbert spaces, the previous theorem can be
reformulated in terms of generalized expanders. Let $\M$ be a class
of metric spaces.

\begin{defn}
A {\it sequence of $(\M,p)$-valued generalized expanders} is a
sequence of finite metric spaces $(X_n)$ satisfying the following
property. For every function $\rho_+:\R_+\to \R_+$, there exist
$K>0,$ and  a sequence $r_n\to \infty$ such that each $X_n$ carries
a symmetric probability measure $\mu_n$ satisfying
\begin{itemize}
\item[-] $\mu_n$ is supported on $\Delta_n(X_n)$;
\item[-] for all maps $F$ from $X_n$ to a metric space $Y\in \M$, satisfying $$|F(x)-F(y)|\leq \rho_+(d(x,y))\quad \forall(x,y)\in
\Delta_1(X_n),$$ we have
    \begin{equation}\label{expander}
\sum_{x,y\in X_n}|F(x)-F(y)|^p\mu_n(x,y)\leq K^p.
    \end{equation}
\end{itemize}
\end{defn}

\begin{cor}\label{Maincor}
Let $\M$ be a class of metric spaces such that the corresponding
sheafs are $p$-admissible for some $1\leq p<\infty.$ Then a metric
space $X$ does not coarsely embed into any element of $\M$ if and
only if it has a coarsely embedded sequence of $(\M,p)$-valued
generalized expanders. \epr
\end{cor}

\subsection{The invariant setting}

If $G$ is a countable discrete group, a length function on $G$ is a
function $L:G\to \R_+$ satisfying $L(1)=0$, $L(gh)\leq L(g)+L(h)$,
and $L(g^{-1})=L(g)$ for all $g,h\in G.$ Clearly, a length function
on $G$ gives rise to a left-invariant metric
$\sigma_L(g,h)=L(g^{-1}h)$. Conversely, given a left-invariant
metric $\sigma$, we define a length function by $L(g)=\sigma(1,g).$

We can define sheaves of length functions as we defined sheaves of
metrics.

\begin{defn}
A {\it sheaf of length functions} $\F$ on a group $G$ is a family of
pairs $(\Omega, L)$, where $\Omega$ is a symmetric neighborhood of
$1$, and $L:\Omega\to \R_+$ satisfying $L(1)=0,$ $L(gh)\leq
L(g)+L(h)$, $L(g^{-1})=L(g)$ for all $g,h\in \Omega$ such that
$gh\in \Omega$.
\end{defn}

Note that a sheaf of lengths naturally induces a sheaf of ``locally
invariant" metrics on $G$ by the relation $\sigma(g,h)=L(g^{-1}h)$
(whenever this is well defined). We will say that a sheaf of lengths
is $p$-admissible if so is the corresponding sheaf of metrics.

\begin{exe}
If $\M$ is a class of metric spaces, the sheaf of $\M$-lengths on
$G$ is the set of $(\Omega,L)$ as above, where $L(g)=\sigma(1,g)$
for some $\M$-metric $\sigma$ defined on $\Omega^2$, satisfying
$\sigma(hg,hg')=\sigma(g,g')$ for all $g,g',h\in G$ such that
$g,g',hg,hg'\in \Omega$.
\end{exe}

Let $G$ be a group equipped with a length function $L_0$. Then, one
sees immediately that the proof of Theorem~\ref{MainThm} can be
formulated with length functions instead of metrics, which yields
the following result

\begin{thm}\label{MainThmEquiv}
Let $\F$ be a $p$-admissible sheaf of length functions on $G$. Then
there exists no coarse length in $\F(G)$ if and only if there exist
$K>0,$ and for every $n\in \N$, a symmetric, finitely supported
probability measure $\mu_n$ on $\{g,L_0(g)\geq n\}$ with the
following property: for every $L\in \F(G)$ satisfying $$L(g)\leq
L_0(g), \quad \forall g\in G,$$ one has
$$\E_{\mu_n} (L^p)\leq K^p.$$
\end{thm}

\subsection{Quantitative statements}
Let $\F$ be a sheaf of metrics on $X$. The $\F$-compression rate of
$X$, denoted by $R_{\F}(X)$ is the supremum of all $\alpha>0$ such
that there exists $\sigma\in \F(X)$ satisfying $d(x,y)^{\alpha}\leq
\sigma(x,y)\leq d(x,y)$, for $d(x,y)$ large enough. The Hilbert
compression rate, usually denoted by $R(X)$, has been introduced in
\cite{GK} and studied by many authors since then as it provides an
interesting quasi-isometry invariant of finitely generated groups.

A slight modification of the proof of Theorem~\ref{MainThm} yields

\begin{thm}\label{quantitativeThm}
Let $X$ be a metric space. Let $\F$ be a $p$-admissible sheaf of
metrics on $X$. The $\F$-compression rate of $X$ is at most $\alpha$
if and only if for all $\beta>\alpha$, there exist $K>0$, and for
every $n\in N$, a symmetric, finitely supported probability measure
$\mu_n$ on $\Delta_n(X)$, with the following property: for every
$\sigma\in \F(X)$ satisfying $$\sigma(x,y)\leq d(x,y), \quad \forall
(x,y)\in \Delta_1(X),$$ one has
$$\E_{\mu_n} (\sigma^p)\leq (Kn^{\beta})^p.$$
\end{thm}

Assume that $X=G$ is a finitely generated group equipped with a word
metric, denoted by $|g|=|g|_S$, associated to a finite symmetric
generating subset $S$. Theorem~\ref{quantitativeThm} becomes

\begin{thm}\label{quantitativeThmEquiv}
Let $\F$ be a $p$-admissible sheaf of length functions on $G$. The
$\F$-compression rate of $G$ is at most $\alpha$ if and only if
there exist $K>0,$ and for every $n\in \N$, a  symmetric, finitely
supported probability measure $\mu_n$ on $\{g,|g|\geq n\}$ with the
following property: for every $L\in \F(G)$ satisfying $$L(g)\leq
|g|, \quad \forall g\in G$$ one has
$$\E_{\mu_n} (L^p)\leq (K|g|^{\beta})^p.$$
\end{thm}

\section{Proof of Theorem~\ref{MainThm}}

The ``if" part is obvious, as the condition (\ref{maineq}) roughly
says that the sequence $\mu_n$ selects pairs $(x,y)$ of arbitrarily
distant points in $X$ for which $\sigma(x,y)\leq K$.

Let $X$ be a metric space, and let $\F$ be a $p$-admissible sheaf of
metrics on $X.$ We assume that $\F(X)$ contains no coarse metric.
Let $\rho_+:\R_+\to \R_+$ be an increasing unbounded function.

Our first step is the next lemma. Let $\A$ be the set of finite
subsets of $X$ containing a distinguished point $o$.

\begin{lem}\label{Mainlemma}
Assume that there exists a function $T:\R_+\to \R_+$ with the
following property: for all $U\in \A$, and all $K>0$, there exists
$\sigma\in\F(U)$ such that
$$\sigma(x,y)\leq \rho_+(d(x,y)),\quad\forall x,y\in U$$ and
$$\sigma(x,y)\geq T(K)$$ for all $(x,y)\in \Delta_T(U)$.
Then there exists a coarse element in $\F(X)$.
\end{lem}
\bpr As $\F$ is $p$-admissible, up to taking a pointwise limit with
respect to an ultrafilter on $\A$, we can assume that for all $K>0$,
there exists $\sigma \in \F(X)$ such that
$$\sigma(x,y)\leq \rho_+(d(x,y)),\quad\forall x,y\in X$$ and
$$\sigma(x,y)\geq K$$ for all $(x,y)\in \Delta_{T(K)}(X)$.

Let $K_n$ be an increasing sequence satisfying
$$\sum_{n=1}^{\infty}\frac{1}{K_n^p}\leq 1,$$
and let $T$ be as in the lemma. Now take $\sigma_n$ as above, and define
$$\sigma(x,y)=\left(\sum_{n\geq 1}\left(\sigma_n(x,y)/K_n\right)^p\right)^{1/p}$$ for all $x,y\in X$.
The fact that $\sigma$ is well defined follows from the fact that
$\F$ is $p$-admissible. Moreover, we have
$$\rho_-(d(x,y))\leq \sigma(x,y)\leq \rho_+(d(x,y)),$$
for all $x,y\in X,$ where $$\rho_-(t)=\card\{n, T(K_n)\leq t\}.$$
Clearly, $\rho_-(t)\to\infty$ when $t\to\infty$, so we are done.
\epr

\

The second step of the proof is an adaptation of the proof of
\cite[Proposition~15.5.2]{M}. Suppose that $X$ does not coarsely
embed into a Hilbert space.

By the lemma, there exists a number $K$ with the following property.
For all $T$, there exist $U\in A$ such that for all $\sigma\in
\F(U)$ satisfying
\begin{equation}\label{C_2-cond}
\sigma(x,y)\leq \rho_+(d(x,y)),
\end{equation}
there are two points $x,y$ in $U$ such that $T\leq d(x,y)$ and
$$\sigma(d(x,y))<K.$$
Note that we can take $T$ such that $\rho_+(T)\geq K.$

Consider the two following convex subsets of $\ell^2\left(U^2\right).$
Let $C_1$ be the set of functions $\phi:U^2\to \R_+$ satisfying
$$\phi(x,y)\leq \rho_+(d(x,y))^p,\quad \forall (x,y)\in U^2,$$
and
\begin{equation}\label{C_1-cond}
\phi(x,y)\geq K^p, \quad \forall (x,y)\in\Delta_{T}(U).
\end{equation}
Let $C_2$ be the set of $\sigma^p$, where $\sigma\in \F(U)$ satisfies
$$\sigma(x,y)\leq \rho_+(d(x,y)).$$
The previous reformulation of the lemma implies that these two
convex subsets are disjoint. We have even better. For every subset
$V$ of a vector space $E$, we denote
$$\R_+V=\{tv,\: t\in \R_+, v\in V\}.$$
\begin{lem}
The cones $\R_+C_1$ and $\R_+C_2$ intersect only at $\{0\}$.
Moreover, $\{0\}$ is extremal in both cones.
\end{lem}
\bpr  The fact that $\{0\}$ is extremal is just a consequence of the
fact that the two cones only contain non-negative functions. Let
$t>0$ and let $\phi\in C_2\setminus\{0\}$. We want to prove that
$t\phi$ does not belong to $C_1$. By the first condition of
$p$-admissibility, there exists $\sigma_t\in \F(U)$ such that
$t\phi=\sigma_t^p$. Moreover, if $t\phi$ also satisfies
(\ref{C_1-cond}), then $\sigma_t$ satisfies (\ref{C_2-cond}), so
$t\phi\in C_2$, and hence it cannot be in $C_1.$ \epr

\

Hence by Hahn-Banach's theorem, there exists a  vector $u\in \ell^2(U^2)$ such that
$$\langle \phi,u \rangle> 0,$$ for all non-zero $\phi\in C_1$ and
$$\langle \phi,u \rangle < 0,$$ for all non-zero $\phi\in C_2$.

Let $u_+=\max\{u,0\}$ and $u_-=\min\{u,0\}$. One sees that $u_+$ is
non-zero in restriction to $\Delta_T(U)$ by applying the first
inequality to the function
$$\phi(x,y)=\left\{
  \begin{array}{ll}
    0, & \hbox{if $d(x,y)<T$;} \\

    K^p, & \hbox{otherwise.}
  \end{array}
\right.$$

Now, apply the first inequality to the function
$$\phi(x,y)=\left\{
  \begin{array}{ll}
    K^p, & \hbox{if $u(x,y)>0$, and $d(x,y)\geq T$;} \\
    0, & \hbox{if $u(x,y)>0$, and $d(x,y)<T$}; \\
    \rho_+(d(x,y))^p, & \hbox{otherwise.}
  \end{array}
\right.$$ Note that this is possible as $T$ has been chosen such
that $\rho_+(T)\geq K$. We get
$$\sum_{x,y}\rho_+(d(x,y))^pu_-(x,y)\leq K^p\sum_{(x,y)\in \Delta_T(U)}u_+(x,y).$$
On the other hand, if $\phi\in C_2$, i.e.  $\phi(x,y)=\sigma(x,y)^p$, then using the second inequality,
\begin{eqnarray*}
\sum_{(x,y)\in \Delta_T(B(o,T))}\sigma(x,y)^pu_+(x,y) &\leq & \sum_{(x,y)\in B(o,T)^2}\sigma(x,y)^pu_+(x,y)\\
&\leq & \sum_{x,y}\sigma(x,y)^pu_-(x,y)\\
                                & \leq & \sum_{x,y}\rho_+(d(x,y))^pu_-(x,y).
\end{eqnarray*}
Now, combining these two inequalities, we get
$$\sum_{x,y}\sigma(x,y)^p\frac{u_+(x,y)}{\sum_{(x,y)\in \Delta_T(U)}u_+(x,y)}\leq K^p.$$
So the theorem follows by taking the probability measure on $\Delta_T(U)$ defined by $$\mu(x,y)=\frac{u_+(x,y)}{\sum_{(x,y)\in \Delta_T(U)}u_+(x,y)}.\;\blacksquare$$

\bibliographystyle{amsplain}

\bigskip

\footnotesize

\noindent \noindent Romain Tessera\\
Department of mathematics, Vanderbilt University,\\ Stevenson
Center, Nashville, TN 37240. \\ E-mail:
\url{romain.a.tessera@vanderbilt.edu}

\end{document}